\documentclass[12pt]{article}
\usepackage{amssymb,amsmath,amsthm,amsfonts,enumerate, bbm, mathdots}

\usepackage{latexsym}
\usepackage{amscd}
\usepackage{stmaryrd}
\usepackage{color}
\usepackage[all]{xypic}
\usepackage{epsfig}
\usepackage{graphics}
\usepackage{ifthen}
\usepackage{varioref}
\usepackage{rotating}
\usepackage{extarrows}
\usepackage{cite}
\usepackage{mathrsfs}

\numberwithin{equation}{section}

\theoremstyle{plain}
\newtheorem{thm}{Theorem}[section]
\newtheorem{cor}[thm]{Corollary}
\newtheorem{lem}[thm]{Lemma}
\newtheorem{prop}[thm]{Proposition}

\definecolor{darkgreen}{rgb}{0.0625,0.64,0.0625}
\usepackage[
            pdfstartview=FitH,
            CJKbookmarks=true,
            bookmarksnumbered=true,
            bookmarksopen=true,
            colorlinks,
            pdfborder=001,
            linkcolor=blue,
            anchorcolor=green,
            citecolor=red
            ]{hyperref}

\newfont{\scyr}{wncyr10 scaled 550}
\def\shuffle{\,\mbox{\bf \scyr X}\,}

\def\proof{\noindent {\bf Proof.\;}}

\def\reg{\operatorname{reg}}

\allowdisplaybreaks

\begin{document}

\title{Weighted sum formula of multiple $L$-values and its applications}

\date{\small ~ \qquad\qquad School of Mathematical Sciences, Tongji University \newline No. 1239 Siping Road,
Shanghai 200092, China}

\author{Zhonghua Li\thanks{E-mail address: zhonghua\_li@tongji.edu.cn} ~and ~Zhenlu Wang\thanks{E-mail address: zlw@tongji.edu.cn}}

\maketitle

\begin{abstract}
In this paper, we study the multiple $L$-values and the multiple zeta values of level $N$. We set up the algebraic framework for the double shuffle relations of the multiple zeta values of level $N$. Using the regularized double shuffle relations of multiple $L$-values, we give a sum formula and a weighted sum formula of multiple $L$-values. As applications, we give sum formulas and weighted sum formulas of double zeta values of level $2$ and $3$.
\end{abstract}

{\small

}


\section{Introduction}\label{Sec:Intro}

In this paper, we study the multiple $L$-values and the multiple zeta values of level $N$.
We recall the definition of multiple $L$-values from \cite{Arakawa-Kaneko}. Let $N$ be a fixed positive integer. Set $R=R_N=\mathbb{Z}/N\mathbb{Z}$. Fix a primitive $N$th root of unity $\omega=\omega_N=\exp(2\pi i/N)$. For positive integers $n,k_1,k_2,\ldots,k_n$ and $a_1,a_2,\ldots,a_n\in R$, the multiple $L$-values are defined by
\begin{align*}
&L_{\ast}(k_1,\ldots,k_n;a_1,\ldots,a_n)=\sum\limits_{m_1>\cdots>m_n>0}\frac{\omega^{a_1m_1+\cdots+a_nm_n}}{m_1^{k_1}\cdots m_n^{k_n}},\\
&L_{\shuffle}(k_1,\ldots,k_n;a_1,\ldots,a_n)=\sum\limits_{m_1>\cdots>m_n>0}\frac{\omega^{a_1(m_1-m_2)+\cdots+a_{n-1}(m_{n-1}-m_n)+a_nm_n}}{m_1^{k_1}\cdots m_n^{k_n}}.
\end{align*}
The above series are convergent when $k_1\geq 2$ or $k_1=1$ and $a_1\neq 0$. We have the relations
\begin{align*}
&L_{\ast}(k_1,\ldots,k_n;a_1,\ldots,a_n)=L_{\shuffle}(k_1,\ldots,k_n;a_1,a_1+a_2,\ldots,a_{1}+\cdots+a_n),\\
&L_{\shuffle}(k_1,\ldots,k_n;a_1,\ldots,a_n)=L_{\ast}(k_1,\ldots,k_n;a_1,a_2-a_1,\ldots,a_n-a_{n-1}).
\end{align*}
We also have the iterated integral representation
$$L_{\shuffle}(k_1,\ldots,k_n;a_1,\ldots,a_n)=\int_0^1\left(\frac{dt}{t}\right)^{k_1-1}\frac{\omega^{a_1}dt}{1-\omega^{a_1}t}\cdots\left(\frac{dt}{t}\right)^{k_n-1}\frac{\omega^{a_n}dt}{1-\omega^{a_n}t},$$
where for one forms $\omega_i=f_i(t)dt$, we define
$$\int_0^1\omega_1\omega_2\cdots\omega_k=\int\limits_{1>t_1>t_2>\cdots>t_k>0}f_1(t_1)f_2(t_2)\cdots f_k(t_k)dt_1dt_2\cdots dt_k.$$

If $N=1$, the multiple $L$-values are just multiple zeta values
$$\zeta(k_1,\ldots,k_n)=\sum\limits_{m_1>\cdots>m_n>0}\frac{1}{m_1^{k_1}\cdots m_n^{k_n}},$$ where $k_1,\ldots,k_n\in \mathbb{N}$ with $k_1\geq 2$. If $N=2$, then $R=\{0,1\}$ and $\omega=-1$. Therefore the multiple $L$-values are just the alternating multiple zeta values:
$$L_{\ast}(k_1,\ldots,k_n;a_1,\ldots,a_n)=\sum\limits_{m_1>\cdots>m_n>0}\frac{(-1)^{a_1m_1+\cdots+a_nm_n}}{m_1^{k_1}\cdots m_n^{k_n}}.$$
As usual, the above value is denoted by $\zeta^{(2)}\left(\fbox{$k_1$},\ldots,\fbox{$k_n$}\right)$, where
$$\fbox{$k_i$}=\begin{cases}
k_i & \text{if\;} a_i=0,\\
 \overline{k_i} & \text{if\;} a_i=1.
\end{cases}$$
Hence for example, we set $L_{\ast}(k_1,k_2;0,1)=\zeta^{(2)}(k_1,\overline{k_2})$. 
If $N=3$, then $R=\{0,1,2\}$. In this time, we denote $L_{\ast}(k_1,\ldots,k_n;a_1,\ldots,a_n)$ by $\zeta^{(3)}\left(\fbox{$k_1$},\ldots,\fbox{$k_n$}\right)$ with
$$\fbox{$k_i$}=\begin{cases}
k_i & \text{if\;} a_i=0,\\
 \overline{k_i} & \text{if\;} a_i=1,\\
\widetilde{k_i} & \text{if\;} a_i=2.
\end{cases}$$
Hence for example, we set $L_{\ast}(k_1,k_2,k_3;0,1,2)=\zeta^{(3)}(k_1,\overline{k_2},\widetilde{k_3})$. Note that
$$\zeta^{(N)}(k_1,\ldots,k_n)=\zeta(k_1,\ldots,k_n)$$
for $N=2,3$.

In \cite{Guo-Xie2009}, Guo and Xie gave a weighted sum formula of multiple zeta values by the regularized double shuffle relations. By using the regularized double shuffle relations of multiple $L$-values, we obtain a sum formula and a weighted sum formula of multiple $L$-values in this paper.

In \cite{Xu-Zhao}, Xu and Zhao studied a variant of multiple zeta values of level $2$ (which is called multiple mixed values therein), which forms a subspace of the space of alternating multiple zeta values. This variant includes both Hoffman's multiple $t$-values \cite{Hoffman2019} and Kaneko-Tsumura's multiple $T$-values \cite{Kaneko-Tsumura} as special cases. The multiple zeta values of any level is introduced by Yuan and Zhao in \cite{Yuan-Zhao}. For positive integers $n,k_1,k_2,\ldots,k_n$ with $k_1\geq 2$ and $a_1,a_2,\ldots,a_n\in R$, the multiple zeta values of level $N$ is defined by
\begin{align*}
&\zeta_{N}(k_1,\ldots,k_n;a_1,\ldots,a_n)=\sum\limits_{m_1>\cdots>m_n>0\atop m_i\equiv a_i \pmod N}\frac{N^n}{m_1^{k_1}\cdots m_n^{k_n}}.
\end{align*}
In this paper, we set up the algebra framework for the double shuffle relations of multiple zeta values of level $N$.

In Section \ref{Sec:MPL}, we recall the double shuffle relations of multiple $L$-values. In Section \ref{Sec:Level N Mzvs}, we study the double shuffle relations of multiple zeta values of level $N$. In Section \ref{Sec:Formulas}, we give a sum formula and a weighted sum formula of multiple $L$-values. As applications, we provide some sum formulas and weighted sum formulas of double zeta values of level $2$ and level $3$. In particular, we reprove the weighted sum formula of double $T$-values appeared in \cite{Kaneko-Tsumura}.


\section{Double shuffle relations of multiple $L$-values}\label{Sec:MPL}

We recall the double shuffle relations of multiple $L$-values from \cite{Arakawa-Kaneko}. Let $\mathcal{A}=\mathbb{Q}\langle x,y_a\mid a\in R_N\rangle$ be the non-commutative polynomial algebra generated by the alphabet $\{x,y_a\mid a\in R_N\}$. Define the subalgebras
$$\mathcal{A}^1=\mathbb{Q}+\sum\limits_{a\in R_N}\mathcal{A}y_a$$
and
$$\mathcal{A}^0=\mathbb{Q}+\sum\limits_{a\in R_N}x\mathcal{A}y_a+\sum\limits_{a,b\in R_N,b\neq 0}y_b\mathcal{A}y_a.$$
For any $k\in\mathbb{N}$ and $a\in R$, set $z_{k,a}=x^{k-1}y_a$. Define the evaluation maps $\mathcal{L}_{\ast}:\mathcal{A}^0\longrightarrow \mathbb{C}$ and $\mathcal{L}_{\shuffle}:\mathcal{A}^0\longrightarrow \mathbb{C}$ by $\mathbb{Q}$-linearities and
\begin{align*}
&\mathcal{L}_{\ast}(z_{k_1,a_1}\cdots z_{k_n,a_n})=L_{\ast}(k_1,\ldots,k_n;a_1,\ldots,a_n),\\
&\mathcal{L}_{\shuffle}(z_{k_1,a_1}\cdots z_{k_n,a_n})=L_{\shuffle}(k_1,\ldots,k_n;a_1,\ldots,a_n).
\end{align*}
Let $\mathcal{I}:\mathcal{A}^1\longrightarrow\mathcal{A}^1$ be the $\mathbb{Q}$-linear endomorphism of $\mathcal{A}^1$ determined by
$$\mathcal{I}(z_{k_1,a_1}z_{k_2,a_2}\cdots z_{k_n,a_n})=z_{k_1,a_1}z_{k_2,a_1+a_2}\cdots z_{k_n,a_1+\cdots+a_n}.$$
The linear map $\mathcal{I}$ is invertible, and the inverse $\mathcal{I}^{-1}:\mathcal{A}^1\longrightarrow\mathcal{A}^1$ satisfies
$$\mathcal{I}^{-1}(z_{k_1,a_1}z_{k_2,a_2}\cdots z_{k_n,a_n})=z_{k_1,a_1}z_{k_2,a_2-a_1}\cdots z_{k_n,a_n-a_{n-1}}.$$
Hence we have $\mathcal{L}_{\ast}=\mathcal{L}_{\shuffle}\circ \mathcal{I}$, or equivalently  $\mathcal{L}_{\shuffle}=\mathcal{L}_{\ast}\circ \mathcal{I}^{-1}$. More precisely, for any $w\in\mathcal{A}^0$, we have
$$\mathcal{L}_{\ast}(w)=\mathcal{L}_{\shuffle}\left(\mathcal{I}(w)\right)\quad \text{and} \quad \mathcal{L}_{\shuffle}(w)=\mathcal{L}_{\ast}\left(\mathcal{I}^{-1}(w)\right).$$

The harmonic shuffle product $\ast$ on $\mathcal{A}^1$ is defined by $\mathbb{Q}$-bilinearity and the rules
\begin{align*}
&1\ast w=w\ast 1=w,\\
&z_{k,a}w_1\ast z_{l,b}w_2=z_{k,a}(w_1\ast z_{l,b}w_2)+z_{l,b}(z_{k,a}w_1\ast w_2)+z_{k+l,a+b}(w_1\ast w_2),
\end{align*}
for all $k,l\geq 1$, $a,b\in R_N$, and any words $w,w_1,w_2\in\mathcal{A}^1$. The harmonic shuffle product $\ast$ is associative and commutative. Hence we get a commutative algebra $(\mathcal{A}^1,\ast)$ and its subalgebra $(\mathcal{A}^0,\ast)$, which are denoted by $\mathcal{A}^1_{\ast}$ and $\mathcal{A}^0_{\ast}$, respectively. The shuffle product $\shuffle$ on $\mathcal{A}$ is defined by $\mathbb{Q}$-bilinearity and the rules
\begin{align*}
&1\shuffle w=w\shuffle 1=w,\\
&uw_1\shuffle vw_2=u(w_1\shuffle vw_2)+v(uw_1\shuffle w_2),
\end{align*}
for any words $w,w_1,w_2\in\mathcal{A}$ and $u,v\in\{x,y_a\mid a\in R_N\}$. Then we have the commutative algebra $\mathcal{A}_{\shuffle}$ and its subalgebras $\mathcal{A}^1_{\shuffle}$ and $\mathcal{A}^0_{\shuffle}$. For any $w_1,w_2\in\mathcal{A}^0$, we have
$$\mathcal{L}_{\ast}(w_1\ast w_2)=\mathcal{L}_{\ast}(w_1)\mathcal{L}_{\ast}(w_2)\quad \text{and} \quad \mathcal{L}_{\shuffle}(w_1\ast w_2)=\mathcal{L}_{\shuffle}(w_1)\mathcal{L}_{\shuffle}(w_2),$$
which induce the finite double shuffle relations
$$\mathcal{L}_{\shuffle}\left(\mathcal{I}(w_1)\shuffle\mathcal{I}(w_2)-\mathcal{I}(w_1\ast w_2)\right)=0,\quad (w_1,w_2\in\mathcal{A}^0),$$
or equivalently
 $$\mathcal{L}_{\ast}\left(\mathcal{I}^{-1}(w_1)\ast\mathcal{I}^{-1}(w_2)-\mathcal{I}^{-1}(w_1\shuffle w_2)\right)=0,\quad (w_1,w_2\in\mathcal{A}^0).$$

Since $\mathcal{A}^1_{\ast}=\mathcal{A}^0_{\ast}[y_0]$ and $\mathcal{A}^1_{\shuffle}=\mathcal{A}^0_{\shuffle}[y_0]$, one can define the regularization maps $\reg_{\ast}:\mathcal{A}^1_{\ast}\longrightarrow\mathcal{A}^0_{\ast}$ and $\reg_{\shuffle}:\mathcal{A}^1_{\shuffle}\longrightarrow\mathcal{A}^0_{\shuffle}$, which are algebraic morphisms. Hence we have the regularized double shuffle relations
$$\mathcal{L}_{\shuffle}\left(\reg_{\shuffle}\left(\mathcal{I}(w_0\ast w_1)-\mathcal{I}(w_0)\shuffle\mathcal{I}(w_1)\right)\right)=0,\quad (w_0\in\mathcal{A}^0,w_1\in\mathcal{A}^1)$$
and
$$\mathcal{L}_{\ast}\left(\reg_{\ast}\left(\mathcal{I}^{-1}(w_0\shuffle w_1)-\mathcal{I}^{-1}(w_0)\ast\mathcal{I}^{-1}(w_1)\right)\right)=0,\quad (w_0\in\mathcal{A}^0,w_1\in\mathcal{A}^1).$$


\section{Double shuffle relations of multiple zeta values of level $N$}\label{Sec:Level N Mzvs}

\subsection{Multiple zeta values of level $N$}

Recall that $\omega=\omega_N=\exp(2\pi i/N)$ is a fixed primitive $N$th root of unity. The following lemma indicates that the multiple zeta values of level $N$ can be expressed by multiple $L$-values.

\begin{lem}\label{Lem:MZVN-MLV}
For positive integers $n,k_1,k_2,\ldots,k_n$ with $k_1\geq 2$ and $a_1,a_2,\ldots,a_n\in R$, we have
\begin{align*}
&\zeta_{N}(k_1,\ldots,k_n;a_1,\ldots,a_n)\\
=&\sum\limits_{m_1>m_2>\cdots>m_n>0}\frac{\prod_{i=1}^{n}(1+\omega^{m_i-a_i}+\omega^{2(m_i-a_i)}+\cdots+\omega^{(N-1)(m_i-a_i)})}{m_1^{k_1}m_2^{k_2}\cdots m_n^{k_n}}
\end{align*}
\end{lem}

\proof
As
$$1+\omega^r+\omega^{2r}+\cdots+\omega^{(N-1)r}=
\begin{cases}
N & \text{if\;} N\mid r, \\
0 & \text{if\;} N\nmid r,
\end{cases}$$
we get the result.
\qed

Using the series representations, we get the harmonic shuffle structure among the multiple zeta values of level $N$ as displaying in the following simple example:
\begin{align*}
&\zeta_N(k;a)\zeta_N(l;b)=\sum\limits_{m=1 \atop m\equiv a\pmod N}^{\infty}\sum\limits_{n=1 \atop n\equiv b\pmod N}^{\infty}\frac{N^2}{m^{k}n^{l}}\\
&=\sum\limits_{m>n>0 \atop m\equiv a, n\equiv b\pmod N}\frac{N^2}{m^{k}n^{l}}+\sum\limits_{n>m>0 \atop m\equiv a, n\equiv b\pmod N}\frac{N^2}{n^{l}m^{k}}+\delta_{a,b}\sum\limits_{m=1\atop m\equiv a\pmod N}^{\infty}\frac{N^2}{m^{k+l}}\\
&=\zeta_N(k,l;a,b)+\zeta_N(l,k;b,a)+\delta_{a,b}N\zeta_N(k+l;a),
\end{align*}
where $k,l\geq 2$, $a,b\in R$ and $\delta_{a,b}$ is the Kronecker symbol.

To study the shuffle structure among the multiple zeta values of level $N$, we introduce a map $r:\mathbb{Z}\longrightarrow \{1,2,\ldots,N\}$, which is defined  by
$$r(a)\equiv a \pmod N \quad \text{and}\quad r(a)\in\{1,2,\ldots,N\}$$
for any $a\in \mathbb{Z}$. We also define one forms
$$\Omega_0=\frac{dt}{t},\quad \Omega_a=\frac{Nt^{a-1}dt}{1-t^N},$$
where $a\in\{1,2,\ldots,N\}$. Then we have the iterated integral representation.

\begin{lem}
Let $k_1,k_2,\ldots,k_n$ be positive integers with $k_1\geq 2$.
\begin{description}
\item [(1)] For any $b_1,\ldots,b_n\in \{1,2,\ldots,N\}$, we have
\begin{align*}
&\int_{0}^{1}\Omega_0^{k_1-1}\Omega_{b_1}\cdots\Omega_0^{k_n-1}\Omega_{b_n}=\zeta_N(k_1,\ldots,k_n;b_1+\cdots+b_n,\ldots,b_{n-1}+b_n,b_n).
\end{align*}
\item [(2)] For any $a_1,a_2,\cdots,a_n\in R$, we have
\begin{align*}
&\zeta_N(k_1,\ldots,k_n;a_1,\ldots,a_n)=\int_{0}^{1}\Omega_0^{k_1-1}\Omega_{r(a_1-a_2)}\Omega_0^{k_2-1}\Omega_{r(a_2-a_3)}
\cdots\\
&\qquad\qquad\qquad\qquad\qquad\qquad\qquad\times\Omega_0^{k_{n-1}-1}\Omega_{r(a_{n-1}-a_n)}\Omega_0^{k_n-1}\Omega_{r(a_n)}.
\end{align*}
\end{description}
\end{lem}

\proof
We prove (1). As
$$\int_{0}^{t}\frac{t^{b_n-1}dt}{1-t^N}=\sum\limits_{l=0}^{\infty}\int_{0}^{t}t^{lN+b_n-1}dt=\sum\limits_{l=0}^{\infty}\frac{t^{lN+b_n}}{lN+b_n},$$
we get
$$\int_{0}^{t}\left(\frac{dt}{t}\right)^{k_n-1}\frac{t^{b_n-1}dt}{1-t^N}=\sum\limits_{l=0}^{\infty}\frac{t^{lN+b_n}}{(lN+b_n)^{k_n}}.$$
Similarly, as
\begin{align*}
\int_{0}^{t}\frac{t^{b_{n-1}-1}dt}{1-t^N}\left(\frac{dt}{t}\right)^{k_n-1}\frac{t^{b_n-1}dt}{1-t^N}&=\sum\limits_{l_2=0}^{\infty}\frac{1}{(l_2N+b_n)^{k_n}}\sum\limits_{l_1=0}^{\infty}\int_{0}^{t}t^{(l_1+l_2)N+b_{n-1}+b_n-1}dt\\
&=\sum\limits_{l_1,l_2\geq 0}\frac{t^{(l_1+l_2)N+b_{n-1}+b_n}}{(l_2N+b_n)^{k_n}((l_1+l_2)N+b_{n-1}+b_n)},
\end{align*}
we find
\begin{align*}
&\int_{0}^{t}\left(\frac{dt}{t}\right)^{k_{n-1}-1}\frac{t^{b_{n-1}-1}dt}{1-t^N}\left(\frac{dt}{t}\right)^{k_n-1}\frac{t^{b_n-1}dt}{1-t^N}\\
=&\sum\limits_{l_1,l_2\geq 0}\frac{t^{(l_1+l_2)N+b_{n-1}+b_n}}{(l_2N+b_n)^{k_n}((l_1+l_2)N+b_{n-1}+b_n)^{k_{n-1}}}.
\end{align*}
Then by induction, one easily get the result.
\qed

\subsection{Algebraic setup}

Let $\mathcal{U}=\mathbb{Q}\langle x_0,x_1,\ldots,x_N\rangle$ be the non-commutative algebra generated by the alphabet $\{x_a\mid a=0,1,\ldots,N\}$. Define the subalgebras
$$\mathcal{U}^1=\mathbb{Q}+\sum\limits_{a=1}^{N}\mathcal{U}x_a$$ spanned by words not ending with $x_0$
and
$$\mathcal{U}^0=\mathbb{Q}+\sum\limits_{a=1}^{N}x_0\mathcal{U}x_a$$ spanned by words begining with $x_0$ and not ending with $x_0$. We set $y_{k,a}=x_0^{k-1}x_a$, where $k \in \mathbb{N}$ and $a\in \{1,2,\ldots,N\}$.

We define the $\mathbb{Q}$-linear map (called the evaluation map) $\zeta_N:\mathcal{U}^0\longrightarrow \mathbb{R}$ by $\zeta_N(1_x)=1$ and
$$\zeta_N(y_{k_1,a_1}\cdots y_{k_n,a_n})=\zeta_N(k_1,\ldots,k_n;a_1,\ldots,a_n),$$
where $1_x$ is the empty word, $k_1,\ldots,k_n\in \mathbb{N}$, $k_1\geq 2$ and $a_1,\cdots,a_n\in \{1,2,\cdots,N\}$.

We define the stuffle product $\ast$ on $\mathcal{U}^1$ by $\mathbb{Q}$-bilinearity and the rules:
\begin{align*}
&1_x \ast w=w=w \ast 1_x,\\
&y_{k,a}w_1 \ast y_{l,b}w_2=y_{k,a}(w_1\ast y_{l,b}w_2)+y_{l,b}(y_{k,a}w_1\ast w_2)+\delta_{a,b}Ny_{k+l,a}(w_1\ast w_2),
\end{align*}
where $w,w_1,w_2$ are words in $\mathcal{U}^1$, $k,l\in \mathbb{N}$, and $a,b\in \{1,2,\cdots,N\}$. The stuffle product $\ast$ is commutative and associative. Therefore $\mathcal{U}^1$ is a commutative $\mathbb{Q}$-algebra with respect to $\ast$. We denote it by $\mathcal{U}_{\ast}^1$. The subspace $\mathcal{U}^0$ is a subalgebra of $\mathcal{U}^1_{\ast}$ and we denote it by $\mathcal{U}_{\ast}^0$. Then from the infinite series representations of multiple zeta values of level $N$, we have the following result.

\begin{prop}
The map $\zeta_N:\mathcal{U}_{\ast}^0\longrightarrow \mathbb{R}$ is an algebra homomorphism. More precisely, for any $w_1,w_2 \in \mathcal{U}^0$, we have
$$\zeta_N(w_1\ast w_2)=\zeta_N(w_1)\zeta_N(w_2).$$
\end{prop}

The shuffle product $\shuffle$ on $\mathcal{U}$ is defined by $\mathbb{Q}$-bilinearity and the rules
\begin{align*}
&1\shuffle w=w\shuffle 1=w,\\
&uw_1\shuffle vw_2=u(w_1\shuffle vw_2)+v(uw_1\shuffle w_2),
\end{align*}
where $w,w_1,w_2$ are words in $\mathcal{U}$ and $u,v\in\{x_a\mid a=0,1,\ldots,N\}$. Then we have the commutative algebra $\mathcal{U}_{\shuffle}$ and its subalgebras $\mathcal{U}^1_{\shuffle}$ and $\mathcal{U}^0_{\shuffle}$.

Let $\mathcal{J}$ be the $\mathbb{Q}$-linear endomorphism of $\mathcal{U}^1$ determined by
$$\mathcal{J}(y_{k_1,a_1}\cdots y_{k_n,a_n})=y_{k_1,r(a_1-a_2)}y_{k_2,r(a_2-a_3)}\cdots y_{k_{n-1},r(a_{n-1}-a_n)}y_{k_n,r(a_n)},$$
where $k_1,\ldots,k_n \in \mathbb{N}$ and $a_1,\ldots,a_n \in \{1,2,\cdots,N\}$. It is obvious that $\mathcal{J}$ is invertible, and the inverse $\mathcal{J}^{-1}$ satisfies
$$\mathcal{J}^{-1}(y_{k_1,a_1}\cdots y_{k_n,a_n})=y_{k_1,r(a_1+\cdots+a_n)}y_{k_2,r(a_2+\cdots+a_n)}\cdots y_{k_{n-1},r(a_{n-1}+a_n)}y_{k_n,r(a_n)}.$$

Then from the iterated integral representations of multiple zeta values of level $N$, we have the following result.

\begin{prop}
For any $w_1,w_2 \in \mathcal{U}^0$, we have
$$\zeta_N(\mathcal{J}^{-1}(w_1\shuffle w_2))=\zeta_N(\mathcal{J}^{-1}(w_1))\zeta_N(\mathcal{J}^{-1}(w_2)).$$
\end{prop}

Finally, we get the finite double shuffle relations of multiple zeta values of level $N$.

\begin{thm}[\textbf{Finite double shuffle relation}]
For any $w_1,w_2\in \mathcal{U}^0$, we have
$$\zeta_N\left(\mathcal{J}^{-1}(w_1)\ast \mathcal{J}^{-1}(w_2)-\mathcal{J}^{-1}(w_1\shuffle w_2)\right)=0.$$
\end{thm}




\section{Sum formulas and weighted sum formulas}\label{Sec:Formulas}

In this section, using the regularized double shuffle relations, we derive a sum formula and a weighted sum formula of multiple $L$-values. As applications, we (re)obtain some sum formulas and weighted sum formulas of double zeta values of level $2$ and level $3$.

\subsection{Sum and weighted sum formulas of multiple $L$-values}

We first compute the stuffle products.

\begin{lem}\label{Lem:Stuffle-MLV}
For positive integers $k,n$ with $k\geq n+1$, $n\geq 2$, and $a,a_1,\ldots,a_{n-1}\in R$, we have
\begin{align*}
&\sum\limits_{k_1+\cdots+k_{n-1}=k-1\atop k_j\geq 1,k_1\geq 2}\mathcal{I}^{-1}(z_{1,a})\ast\mathcal{I}^{-1}\left(z_{k_1,a_1}\cdots z_{k_{n-1},a_{n-1}}\right)\\
=&\sum\limits_{k_1+\cdots+k_{n-1}=k-1\atop k_j\geq 1,k_1\geq 2}z_{1,a}z_{k_1,a_1}z_{k_2,a_2-a_1}\cdots z_{k_{n-1},a_{n-1}-a_{n-2}}\\
&+\sum\limits_{i=2}^n\sum\limits_{k_1+\cdots+k_n=k\atop k_j\geq 1,k_1\geq 2,k_i=1}z_{k_1,a_1}\cdots z_{k_{i-1},a_{i-1}-a_{i-2}}z_{k_i,a}z_{k_{i+1}.a_i-a_{i-1}}\cdots z_{k_n,a_{n-1}-a_{n-2}}\\
&+\sum\limits_{k_1+\cdots+k_{n-1}=k\atop k_j\geq 1,k_1\geq 3}z_{k_1,a+a_1}z_{k_2,a_2-a_1}\cdots z_{k_{n-1},a_{n-1}-a_{n-2}}\\
&+\sum\limits_{i=2}^{n-1}\sum\limits_{k_1+\cdots+k_{n-1}=k\atop k_j\geq 1,k_1,k_i\geq 2}z_{k_1,a_1}\cdots z_{k_{i-1},a_{i-1}-a_{i-2}}z_{k_i,a+a_i-a_{i-1}}z_{k_{i+1},a_{i+1}-a_i}\cdots z_{k_{n-1},a_{n-1}-a_{n-2}}
\end{align*}
and
\begin{align*}
&\sum\limits_{l+k_1+\cdots+k_{n-1}=k\atop l,k_j\geq 1, k_1\geq 2}\mathcal{I}^{-1}(z_{l,a})\ast\mathcal{I}^{-1}(z_{k_1,a_1}\cdots z_{k_{n-1},a_{n-1}})\\
=&\sum\limits_{k_1+\cdots+k_n=k\atop k_j\geq 1,k_2\geq 2}z_{k_1,a}z_{k_2,a_1}z_{k_3,a_2-a_1}\cdots z_{k_n,a_{n-1}-a_{n-2}}\\
&+\sum\limits_{i=1}^{n-1}\sum\limits_{k_1+\cdots+k_n=k\atop k_j\geq 1,k_1\geq 2}z_{k_1,a_1}z_{k_2,a_2-a_1}\cdots z_{k_i,a_i-a_{i-1}}z_{k_{i+1},a}z_{k_{i+2},a_{i+1}-a_i}\cdots z_{k_n,a_{n-1}-a_{n-2}}\\
&+\sum\limits_{k_1+\cdots+k_{n-1}=k\atop k_j\geq 1,k_1\geq 2}(k_1-2)z_{k_1,a+a_1}z_{k_2,a_2-a_1}\cdots z_{k_{n-1},a_{n-1}-a_{n-2}}\\
&+\sum\limits_{i=2}^{n-1}\sum\limits_{k_1+\cdots+k_{n-1}=k\atop k_j\geq 1,k_1\geq 2}(k_i-1)z_{k_1,a_1}z_{k_2,a_2-a_1}\cdots z_{k_{i-1},a_{i-1}-a_{i-2}}z_{k_i,a+a_i-a_{i-1}}\\
&\qquad\qquad\qquad\qquad \times z_{k_{i+1},a_{i+1}-a_i}\cdots z_{k_{n-1},a_{n-1}-a_{n-2}}.
\end{align*}
\end{lem}

\proof
As
\begin{align*}
&\mathcal{I}^{-1}(z_{l,a})\ast\mathcal{I}^{-1}(z_{k_1,a_1}\cdots z_{k_{n-1},a_{n-1}})=z_{l,a}\ast z_{k_1,a_1}z_{k_2,a_2-a_1}\cdots z_{k_{n-1},a_{n-1}-a_{n-2}}\\
=&\sum\limits_{i=0}^{n-1} z_{k_1,a_1}z_{k_2,a_2-a_1}\cdots z_{k_i,a_i-a_{i-1}}z_{l,a}z_{k_{i+1},a_{i+1}-a_i}\cdots z_{k_{n-1},a_{n-1}-a_{n-2}}\\
&+\sum\limits_{i=1}^{n-1} z_{k_1,a_1}z_{k_2,a_2-a_1}\cdots z_{k_{i-1},a_{i-1}-a_{i-2}}z_{l+k_i,a+a_i-a_{i-1}}z_{k_{i+1},a_{i+1}-a_i}\cdots z_{k_{n-1},a_{n-1}-a_{n-2}},
\end{align*}
we get the result.
\qed

For shuffle products, we have

\begin{lem}\label{Lem:Shuffle-MLV}
For positive integers $k,n$ with $k\geq n+1$ and $n\geq 2$, $a,a_1,\ldots,a_{n-1}\in R$, we have
\begin{align*}
&\sum\limits_{k_1+\cdots+k_{n-1}=k-1\atop k_j\geq 1,k_1\geq 2}z_{1,a}\shuffle z_{k_1,a_1}\cdots z_{k_{n-1},a_{n-1}}\\
=&\sum\limits_{k_1+\cdots+k_n=k\atop k_j\geq 1,k_1+k_2\geq 3}z_{k_1,a}z_{k_2,a_1}\cdots z_{k_n,a_{n-1}}+\sum\limits_{k_1+\cdots+k_{n-1}=k-1\atop k_j\geq 1,k_1\geq 2}z_{k_1,a_1}\cdots z_{k_{n-1},a_{n-1}}z_{1,a}\\
&+\sum\limits_{i=2}^{n-1}\sum\limits_{k_1+\cdots+k_n=k\atop k_j\geq 1,k_1\geq 2}z_{k_1,a_1}z_{k_2,a_2}\cdots z_{k_{i-1},a_{i-1}}z_{k_i,a}z_{k_{i+1},a_i}\cdots z_{k_n,a_{n-1}}
\end{align*}
and
\begin{align*}
&\sum\limits_{l+k_1+\cdots+k_{n-1}=k\atop l,k_j\geq 1, k_1\geq 2}z_{l,a}\shuffle z_{k_1,a_1}\cdots z_{k_{n-1},a_{n-1}}\\
=&\sum\limits_{k_1+\cdots+k_n=k\atop k_j\geq 1,k_2\geq 2}2^{k_1-1}z_{k_1,a}z_{k_2,a_1}\cdots z_{k_n,a_{n-1}}+\sum\limits_{k_1+\cdots+k_n=k\atop k_j\geq 1,k_2=1}(2^{k_1-1}-1)z_{k_1,a}z_{k_2,a_1}\cdots z_{k_n,a_{n-1}}\\
&+\sum\limits_{i=2}^{n-1}\sum\limits_{k_1+\cdots+k_n=k\atop k_j\geq 1}(2^{k_1+\cdots+k_i-i}-2^{k_2+\cdots+k_i-(i-1)})z_{k_1,a_1}\cdots z_{k_{i-1},a_{i-1}}z_{k_{i},a}z_{k_{i+1},a_{i}}\cdots z_{k_n,a_{n-1}}\\
&+\sum\limits_{k_1+\cdots+k_{n}=k\atop k_j\geq 1}(2^{k_1+\cdots+k_{n-1}-(n-1)}-2^{k_2+\cdots+k_{n-1}-(n-2)})z_{k_1,a_1}\cdots z_{k_{n-1},a_{n-1}}z_{k_n,a}.
\end{align*}
\end{lem}

\proof
As
\begin{align*}
&z_{1,a}\shuffle z_{k_1,a_1}\cdots z_{k_{n-1},a_{n-1}}\\
=&\sum\limits_{j=1}^{k_1}z_{j,a}z_{k_1+1-j,a_1}z_{k_2,a_2}\cdots z_{k_{n-1},a_{n-1}}\\
&+\sum\limits_{i=2}^{n-1}\sum\limits_{j=1}^{k_i}z_{k_1,a_1}z_{k_2,a_2}\cdots z_{k_{i-1},a_{i-1}}z_{j,a}z_{k_i+1-j,a_i}z_{k_{i+1},a_{i+1}}\cdots z_{k_{n-1},a_{n-1}}\\
&+z_{k_1,a_1}\cdots z_{k_{n-1},a_{n-1}}z_{1,a},
\end{align*}
we get the first equation.

In general, similarly as in \cite{Li-Qin-shuffle}, we have
\begin{align*}
&z_{l,a}\shuffle z_{k_1,a_1}\cdots z_{k_{n-1},a_{n-1}}\\
=&\sum\limits_{i=1}^{n-1}\sum\limits_{\alpha_1+\cdots+\alpha_{i+1}\atop =l+k_1+\cdots+k_i, \alpha_j\geq 1}\prod\limits_{j=1}^{i-1}\binom{\alpha_j-1}{k_j-1}\binom{\alpha_i-1}{k_i-\alpha_{i+1}}z_{\alpha_1,a_1}\cdots z_{\alpha_{i-1},a_{i-1}}z_{\alpha_i,a}z_{\alpha_{i+1},a_i}\\
&\qquad\qquad\qquad\qquad \times z_{k_{i+1},a_{i+1}}\cdots z_{k_{n-1},a_{n-1}}\\
&+\sum\limits_{\alpha_1+\cdots+\alpha_n\atop =l+k_1+\cdots+k_{n-1},\alpha_j\geq 1}\prod\limits_{j=1}^{n-1}\binom{\alpha_j-1}{k_j-1}z_{\alpha_1,a_1}\cdots z_{\alpha_{n-1},a_{n-1}}z_{\alpha_n,a}.
\end{align*}
Hence we get
$$\sum\limits_{l+k_1+\cdots+k_{n-1}=k\atop l,k_j\geq 1, k_1\geq 2}z_{l,a}\shuffle z_{k_1,a_1}\cdots z_{k_{n-1},a_{n-1}}=S_1+S_2+S_3,$$
where
\begin{align*}
S_1=&\sum\limits_{l+k_1+\cdots+k_{n-1}=k\atop l,k_j\geq 1, k_1\geq 2}\sum\limits_{\alpha_1+\alpha_2=l+k_1\atop \alpha_j\geq 1}\binom{\alpha_1-1}{k_1-\alpha_2}z_{\alpha_1,a}z_{\alpha_{2},a_1}z_{k_2,a_2}\cdots z_{k_{n-1},a_{n-1}},\\
S_2=&\sum\limits_{i=2}^{n-1}\sum\limits_{l+k_1+\cdots+k_{n-1}=k\atop l,k_j\geq 1, k_1\geq 2}\sum\limits_{\alpha_1+\cdots+\alpha_{i+1}\atop =l+k_1+\cdots+k_i, \alpha_j\geq 1}\prod\limits_{j=1}^{i-1}\binom{\alpha_j-1}{k_j-1}\binom{\alpha_i-1}{k_i-\alpha_{i+1}}z_{\alpha_1,a_1}\cdots z_{\alpha_{i-1},a_{i-1}}\\
&\qquad\qquad\qquad\qquad \times z_{\alpha_i,a}z_{\alpha_{i+1},a_i}z_{k_{i+1},a_{i+1}}\cdots z_{k_{n-1},a_{n-1}}\\
S_3=&\sum\limits_{l+k_1+\cdots+k_{n-1}=k\atop l,k_j\geq 1, k_1\geq 2}\sum\limits_{\alpha_1+\cdots+\alpha_n\atop =l+k_1+\cdots+k_{n-1},\alpha_j\geq 1}\prod\limits_{j=1}^{n-1}\binom{\alpha_j-1}{k_j-1}z_{\alpha_1,a_1}\cdots z_{\alpha_{n-1},a_{n-1}}z_{\alpha_n,a}.
\end{align*}

For $S_1$, we have
$$S_1=\sum\limits_{\alpha_1+\alpha_2+k_2+\cdots+k_{n-1}=k\atop \alpha_j,k_p\geq 1}\sum\limits_{k_1\geq 2\atop k_1\geq \alpha_2}\binom{\alpha_1-1}{k_1-\alpha_2}z_{\alpha_1,a}z_{\alpha_{2},a_1}z_{k_2,a_2}\cdots z_{k_{n-1},a_{n-1}}.$$
If $\alpha_2=1$, we get
$$\sum\limits_{k_1\geq 2\atop k_1\geq \alpha_2}\binom{\alpha_1-1}{k_1-\alpha_2}=\sum\limits_{k_1\geq 2}\binom{\alpha_1-1}{k_1-1}=2^{\alpha_1-1}-1.$$
While if $\alpha_2\geq 2$, we have
$$\sum\limits_{k_1\geq 2\atop k_1\geq \alpha_2}\binom{\alpha_1-1}{k_1-\alpha_2}=\sum\limits_{k_1\geq \alpha_2}\binom{\alpha_1-1}{k_1-\alpha_2}=2^{\alpha_1-1}.$$
Hence we find
$$S_1=\sum\limits_{k_1+\cdots+k_n=k\atop k_j\geq 1,k_2\geq 2}2^{k_1-1}z_{k_1,a}z_{k_2,a_1}\cdots z_{k_n,a_{n-1}}+\sum\limits_{k_1+\cdots+k_n=k\atop k_j\geq 1,k_2=1}(2^{k_1-1}-1)z_{k_1,a}z_{k_2,a_1}\cdots z_{k_n,a_{n-1}}.$$
For $S_2$, we have
\begin{align*}
S_2=&\sum\limits_{i=2}^{n-1}\sum\limits_{\alpha_1+\cdots+\alpha_{i+1}+k_{i+1}\atop +\cdots+k_{n-1}=k, \alpha_j,k_p\geq 1}\sum\limits_{k_1=2}^{\alpha_1}\binom{\alpha_1-1}{k_1-1}\prod\limits_{j=2}^{i-1}\sum\limits_{k_j=1}^{\alpha_j}\binom{\alpha_j-1}{k_j-1}\sum\limits_{k_i=\alpha_{i+1}}^{\alpha_i+\alpha_{i+1}-1}\binom{\alpha_i-1}{k_i-\alpha_{i+1}} \\
&\qquad\qquad\qquad\qquad \times z_{\alpha_1,a_1}\cdots z_{\alpha_{i-1},a_{i-1}}z_{\alpha_i,a}z_{\alpha_{i+1},a_i}z_{k_{i+1},a_{i+1}}\cdots z_{k_{n-1},a_{n-1}}.
\end{align*}
Since
$$\sum\limits_{k_1=2}^{\alpha_1}\binom{\alpha_1-1}{k_1-1}=2^{\alpha_1-1}-1,\quad \sum\limits_{k_j=1}^{\alpha_j}\binom{\alpha_j-1}{k_j-1}=2^{\alpha_j-1},\quad \sum\limits_{k_i=\alpha_{i+1}}^{\alpha_i+\alpha_{i+1}-1}\binom{\alpha_i-1}{k_i-\alpha_{i+1}}=2^{\alpha_i-1},$$
we find
$$S_2=\sum\limits_{i=2}^{n-1}\sum\limits_{k_1+\cdots+k_n=k\atop k_j\geq 1}(2^{k_1+\cdots+k_i-i}-2^{k_2+\cdots+k_i-(i-1)})z_{k_1,a_1}\cdots z_{k_{i-1},a_{i-1}}z_{k_{i},a}z_{k_{i+1},a_{i}}\cdots z_{k_n,a_{n-1}}.$$
Similarly, for $S_3$, we have
\begin{align*}
S_3=&\sum\limits_{\alpha_1+\cdots+\alpha_n=k \atop\alpha_j\geq 1}\sum\limits_{k_1=2}^{\alpha_1}\binom{\alpha_1-1}{k_1-1}\prod\limits_{j=2}^{n-1}\sum\limits_{k_j=1}^{\alpha_j}\binom{\alpha_j-1}{k_j-1}z_{\alpha_1,a_1}\cdots z_{\alpha_{n-1},a_{n-1}}z_{\alpha_n,a}\\
=&\sum\limits_{k_1+\cdots+k_{n}=k\atop k_j\geq 1}(2^{k_1+\cdots+k_{n-1}-(n-1)}-2^{k_2+\cdots+k_{n-1}-(n-2)})z_{k_1,a_1}\cdots z_{k_{n-1},a_{n-1}}z_{k_n,a}.
\end{align*}
Then we get the desired result.
\qed

From Lemma \ref{Lem:Stuffle-MLV} and Lemma \ref{Lem:Shuffle-MLV}, we get the following sum formula and weighted sum formula.

\begin{thm}\label{Thm:Sum-MLV}
For positive integers $k,n$ with $k\geq n+1$, $n\geq 2$, and $a,a_1,\ldots,a_{n-1}\in R$, we have
\begin{align*}
&\sum\limits_{k_1+\cdots+k_n=k\atop k_j\geq 1,k_1\geq 2}L_{\ast}(k_1,\ldots,k_n;a,a_1-a,a_2-a_1,\ldots,a_{n-1}-a_{n-2})\\
&+\sum\limits_{i=2}^{n-1}\sum\limits_{k_1+\cdots+k_n=k\atop k_j\geq 1,k_1\geq 2}L_{\ast}(k_1,\ldots,k_n;a_1,a_2-a_1,\ldots,a_{i-1}-a_{i-2},\\
&\qquad\qquad\qquad\qquad\qquad\qquad a-a_{i-1},a_i-a,a_{i+1}-a_i,\ldots,a_{n-1}-a_{n-2})\\
&+\sum\limits_{k_1+\cdots+k_n=k\atop k_1\geq 2,k_n=1}L_{\ast}(k_1,\ldots,k_n;a_1,a_2-a_1,\ldots,a_{n-1}-a_{n-2},a-a_{n-1})\\
=&\sum\limits_{k_1+\cdots+k_n=k\atop k_1=1,k_2\geq 2}\mathcal{L}_{\ast}\left(z_{k_1,a}z_{k_2,a_1}z_{k_3,a_2-a_1}\cdots z_{k_n,a_{n-1}-a_{n-2}}-z_{k_1,a}z_{k_2,a_1-a}z_{k_3,a_2-a_1}\cdots z_{k_n,a_{n-1}-a_{n-2}}\right)\\
+&\sum\limits_{i=2}^n\sum\limits_{k_1+\cdots+k_n=k\atop k_j\geq 1,k_1\geq 2,k_i=1}L_{\ast}(k_1,\ldots,k_n;a_1,a_2-a_1,\ldots,a_{i-1}-a_{i-2},a,a_i-a_{i-1},\ldots,a_{n-1}-a_{n-2})\\
&+\sum\limits_{k_1+\cdots+k_{n-1}=k\atop k_j\geq 1,k_1\geq 3}L_{\ast}(k_1,\ldots,k_{n-1};a+a_1,a_2-a_1,\ldots,a_{n-1}-a_{n-2})\\
+&\sum\limits_{i=2}^{n-1}\sum\limits_{k_1+\cdots+k_{n-1}=k\atop k_j\geq 1,k_1,k_i\geq 2}L_{\ast}(k_1,\ldots,k_{n-1};a_1,a_2-a_1,\ldots,a_{i-1}-a_{i-2},\\
&\qquad\qquad\qquad\qquad\qquad\qquad\qquad\qquad a+a_i-a_{i-1},a_{i+1}-a_i,\ldots,a_{n-1}-a_{n-2}).
\end{align*}
\end{thm}

\begin{thm}\label{Thm:DoubleShuffle-MLV}
For positive integers $k,n$ with $k\geq n+1$ and $n\geq 2$, $a,a_1,\ldots,a_{n-1}\in R$, we have
\begin{align*}
&\sum\limits_{k_1+\cdots+k_n=k\atop k_j\geq 1,k_2\geq 2}\mathcal{L}_{\ast}\left(2^{k_1-1}z_{k_1,a}z_{k_2,a_1-a}z_{k_3,a_2-a_1}\cdots z_{k_n,a_{n-1}-a_{n-2}}\right.\\
&\qquad\qquad\qquad\left.-z_{k_1,a}z_{k_2,a_1}z_{k_3,a_2-a_1}\cdots z_{k_n,a_{n-1}-a_{n-2}}\right)\\
&+\sum\limits_{k_1+\cdots+k_n=k\atop k_j\geq 1,k_1\geq 2,k_2=1}(2^{k_1-1}-1)L_{\ast}(k_1,\ldots,k_n;a,a_1-a,a_2-a_1,\ldots,a_{n-1}-a_{n-2})\\
&+\sum\limits_{i=2}^{n-1}\sum\limits_{k_1+\cdots+k_n=k\atop k_j\geq 1,k_1\geq 2}(2^{k_1+\cdots+k_i-i}-2^{k_2+\cdots+k_i-(i-1)})L_{\ast}(k_1,\ldots,k_n;a_1,a_2-a_1,\ldots,a_{i-1}-a_{i-2},\\
&\qquad\qquad\qquad\qquad a-a_{i-1},a_{i}-a,a_{i+1}-a_i,\ldots,a_{n-1}-a_{n-2})\\
&+\sum\limits_{k_1+\cdots+k_{n}=k\atop k_j\geq 1,k_1\geq 2}(2^{k_1+\cdots+k_{n-1}-(n-1)}-2^{k_2+\cdots+k_{n-1}-(n-2)})L_{\ast}(k_1,\ldots,k_n;a_1,a_2-a_1,\\
&\qquad\qquad\qquad\qquad \ldots,a_{n-1}-a_{n-2},a-a_{n-1})\\
&=\sum\limits_{i=1}^{n-1}\sum\limits_{k_1+\cdots+k_n=k\atop k_j\geq 1,k_1\geq 2}L_{\ast}(k_1,\ldots,k_n;a_1,a_2-a_1,\ldots,a_i-a_{i-1},a,a_{i+1}-a_i,\ldots,a_{n-1}-a_{n-2})\\
&+\sum\limits_{k_1+\cdots+k_{n-1}=k\atop k_j\geq 1,k_1\geq 2}(k_1-2)L_{\ast}(k_1,\ldots,k_{n-1};a+a_1,a_2-a_1,\ldots,a_{n-1}-a_{n-2})\\
&+\sum\limits_{i=2}^{n-1}\sum\limits_{k_1+\cdots+k_{n-1}=k\atop k_j\geq 1,k_1\geq 2}(k_i-1)L_{\ast}(k_1,\ldots,k_{n-1};a_1,a_2-a_1,\ldots,a_{i-1}-a_{i-2},a+a_i-a_{i-1},\\
&\qquad\qquad\qquad\qquad a_{i+1}-a_i,\ldots,a_{n-1}-a_{n-2}).
\end{align*}
\end{thm}

Let $n=2$. From Theorem \ref{Thm:Sum-MLV}, we get the following sum formula of double $L$-values.

\begin{cor}\label{Cor:Sum-DLV}
For an integer $k$ with $k\geq 3$, and $a_1,a_2\in R$, we have
\begin{align*}
\sum\limits_{j=2}^{k-1}L_{\ast}(j,k-j;a_1,a_2)=&L_{\ast}(k-1,1;a_1+a_2,a_1)-L_{\ast}(k-1,1;a_1+a_2,-a_2)\\
&+\mathcal{L}_{\ast}(z_{1,a_1}z_{k-1,a_1+a_2}-z_{1,a_1}z_{k-1,a_2})+L_{\ast}(k,2a_1+a_2).
\end{align*}
\end{cor}

And from Theorem \ref{Thm:DoubleShuffle-MLV}, we get the following weighted sum formula of double $L$-values.

\begin{cor}\label{Cor:Weightedsum-double}
For an integer $k$ with $k\geq 3$, and any $a_1,a_2\in R$, we have
\begin{align*}
&\sum\limits_{j=2}^{k-1}\left(2^{j-1}L_{\ast}(j,k-j;a_1,a_2-a_1)+(2^{j-1}-1)L_{\ast}(j,k-j;a_2,a_1-a_2)\right.\\
&\qquad\left.-L_{\ast}(j,k-j;a_1,a_2)-L_{\ast}(j,k-j;a_2,a_1)\right)\\
=&L_{\ast}(k-1,1;a_1,a_2-a_1)-L_{\ast}(k-1,1;a_1,a_2)\\
&+\mathcal{L}_{\ast}(z_{1,a_1}z_{k-1,a_2}-z_{1,a_1}z_{k-1,a_2-a_1})+(k-2)L_{\ast}(k;a_1+a_2).
\end{align*}
\end{cor}

\subsection{Sum and weighted sum formulas of double zeta values of level $2$}

Setting $N=2$ in Corollary \ref{Cor:Sum-DLV}, and taking all possible values of $(a_1,a_2)$, we get the sum formulas of alternating double zeta values.

\begin{cor}\label{Cor:Sum-Double}
For an integer $k$ with $k\geq 3$, we have
\begin{align*}
&\sum\limits_{j=2}^{k-1}\zeta(j,k-j)=\zeta(k),\\
&\sum\limits_{j=2}^{k-1}\zeta^{(2)}(j,\overline{k-j})=\zeta^{(2)}(\overline{k-1},1)-\zeta^{(2)}(\overline{k-1},\overline{1})+\zeta^{(2)}(\overline{k}),\\
&\sum\limits_{j=1}^{k-1}\zeta^{(2)}(\overline{j},\overline{k-j})=\zeta^{(2)}(\overline{1},k-1)+\zeta^{(2)}(\overline{k}),\\
&\sum\limits_{j=1}^{k-1}\zeta^{(2)}(\overline{j},k-j)=\zeta^{(2)}(\overline{k-1},\overline{1})-\zeta^{(2)}(\overline{k-1},1)+\zeta^{(2)}(\overline{1},\overline{k-1})+\zeta(k).
\end{align*}
\end{cor}

Now let $N=2$ in Corollary \ref{Cor:Weightedsum-double}. In the case of $(a_1,a_2)=(0,0)$, we get
\begin{align}
\sum\limits_{j=2}^{k-1}(2^j-3)\zeta(j,k-j)=(k-2)\zeta(k).
\label{Eq:Double-00}
\end{align}
In the case of $(a_1,a_2)=(0,1)$, we have
\begin{align}
&\sum\limits_{j=2}^{k-1}(2^{j-1}-1)\zeta^{(2)}(j,\overline{k-j})+\sum\limits_{j=2}^{k-1}(2^{j-1}-1)\zeta^{(2)}(\overline{j},\overline{k-j})\nonumber\\
&\qquad-\sum\limits_{j=2}^{k-1}\zeta^{(2)}(\overline{j},k-j)=(k-2)\zeta^{(2)}(\overline{k}).
\label{Eq:Double-01}
\end{align}
In the case of $(a_1,a_2)=(1,0)$, we get
\begin{align}
&\sum\limits_{j=1}^{k-1}2^{j-1}\zeta^{(2)}(\overline{j},\overline{k-j})+\sum\limits_{j=2}^{k-1}(2^{j-1}-2)\zeta^{(2)}(j,\overline{k-j})-\sum\limits_{j=1}^{k-2}\zeta^{(2)}(\overline{j},k-j)\nonumber\\
=&\zeta^{(2)}(\overline{k-1},\overline{1})+(k-2)\zeta^{(2)}(\overline{k}).
\label{Eq:Double-10}
\end{align}
In the case of $(a_1,a_2)=(1,1)$, we get
\begin{align}
&\sum\limits_{j=2}^{k-1}(2^{j}-1)\zeta^{(2)}(\overline{j},k-j)-2\sum\limits_{j=2}^{k-1}\zeta^{(2)}(\overline{j},\overline{k-j})\nonumber\\
=&\zeta^{(2)}(\overline{k-1},1)-\zeta^{(2)}(\overline{k-1},\overline{1})+\zeta^{(2)}(\overline{1},\overline{k-1})-\zeta^{(2)}(\overline{1},k-1)+(k-2)\zeta(k).
\label{Eq:Double-11}
\end{align}

Using the sum formulas of alternating double zeta values and \eqref{Eq:Double-00}-\eqref{Eq:Double-11}, we get the following weighted sum formulas of alternating double zeta values.

\begin{cor}
For an integer $k$ with $k\geq 3$, we have
\begin{align}
&\sum\limits_{j=2}^{k-1}2^j\zeta(j,k-j)=(k+1)\zeta(k),
\label{Eq:WeightedSum-00}\\
&\sum\limits_{j=2}^{k-1}2^j\zeta^{(2)}(j,\overline{k-j})+\sum\limits_{j=2}^{k-1}2^j\zeta^{(2)}(\overline{j},\overline{k-j})=2\zeta(k)+2k\zeta^{(2)}(\overline{k}),
\label{Eq:WeightedSum-01-11}\\
&\sum\limits_{j=2}^{k-1}2^j\zeta^{(2)}(\overline{j},k-j)=(k-1)\zeta(k)+2\zeta^{(2)}(\overline{k}).
\label{Eq:WeightedSum-10}
\end{align}
\end{cor}

We know that double zeta values of level $2$ can be represented by the alternating double zeta values as in the following way:
\begin{align}
\begin{pmatrix}
\zeta_2(k,l)\\
\zeta_2(\overline{k},l)\\
\zeta_2(k,\overline{l})\\
\zeta_2(\overline{k},\overline{l})
\end{pmatrix}=\begin{pmatrix}
1 & 1 & 1 & 1\\
1 & -1 & 1 & -1\\
1 & 1 & -1 & -1\\
1 & -1 & -1 & 1
\end{pmatrix}\begin{pmatrix}
\zeta(k,l)\\
\zeta^{(2)}(\overline{k},l)\\
\zeta^{(2)}(k,\overline{l})\\
\zeta^{(2)}(\overline{k},\overline{l})
\end{pmatrix}.
\label{Eq:Double-Level-2}
\end{align}
Then using Corollary \ref{Cor:Sum-Double}, \eqref{Eq:Double-Level-2} and the fact $\zeta^{(2)}(\overline{k})=(2^{1-k}-1)\zeta(k)$, we get the sum formulas of double zeta values of level $2$.

\begin{cor}
For an integer $k$ with $k\geq 3$, we have
\begin{align*}
&\sum\limits_{j=2}^{k-1}\zeta_2(j,k-j)=\frac{1}{2^{k-2}}\zeta(k),\\
&\sum\limits_{j=2}^{k-1}\zeta_2(\overline{j},k-j)=2\left(\zeta^{(2)}(\overline{k-1},1)-\zeta^{(2)}(\overline{k-1},\overline{1})\right)\\
&\sum\limits_{j=2}^{k-1}\zeta_2(j,\overline{k-j})=2\left(\zeta^{(2)}(\overline{k-1},\overline{1})+\zeta^{(2)}(\overline{1},\overline{k-1})-\zeta^{(2)}(\overline{k-1},1)-\zeta^{(2)}(\overline{1},k-1)\right)\\
&\qquad\qquad\qquad\qquad+4\left(1-\frac{1}{2^k}\right)\zeta(k),\\
&\sum\limits_{j=2}^{k-1}\zeta_2(\overline{j},\overline{k-j})=2\left(\zeta^{(2)}(\overline{1},k-1)-\zeta^{(2)}(\overline{1},\overline{k-1})\right).
\end{align*}
\end{cor}

Similarly, using Corollary \ref{Cor:Sum-Double}, \eqref{Eq:Double-Level-2} and \eqref{Eq:WeightedSum-00}-\eqref{Eq:WeightedSum-10}, we get the weighted sum formulas of double zeta values of level $2$.

\begin{cor}
For an integer $k$ with $k\geq 3$, we have
\begin{align*}
&\sum\limits_{j=2}^{k-1}2^j\zeta_2(j,k-j)=\frac{k+1}{2^{k-2}}\zeta(k),\\
&\sum\limits_{j=2}^{k-1}2^j\zeta_2(j,\overline{k-j})=4(k-1)\left(1-\frac{1}{2^k}\right)\zeta(k).
\end{align*}
\end{cor}

\proof We get the first equation by \eqref{Eq:WeightedSum-00}$+$\eqref{Eq:WeightedSum-01-11}$+$\eqref{Eq:WeightedSum-10}, and the second equation by \eqref{Eq:WeightedSum-00}$+$\eqref{Eq:WeightedSum-10}$-$\eqref{Eq:WeightedSum-01-11}.\qed


\subsection{Sum and weighted sum formulas of double zeta values of level $3$}

We consider the condition of $N=3$.

Taking all possible values of $(a_1,a_2)$ in  Corollary \ref{Cor:Sum-DLV}, we get the following sum formulas of multiple $L$-values of level $3$.

\begin{cor}\label{Cor:Sum-Double of level 3}
For an integer $k$ with $k\geq 3$, we have
\begin{align*}
&\sum\limits_{j=2}^{k-1}\zeta(j,k-j)=\zeta(k),\\
&\sum\limits_{j=2}^{k-1}\zeta^{(3)}(\overline{j},k-j)=\zeta^{(3)}(\overline{1},\overline{k-1})-\zeta^{(3)}(\overline{1},k-1)+\zeta^{(3)}(\overline{k-1},\overline{1})-\zeta^{(3)}(\overline{k-1},1)+\zeta^{(3)}(\widetilde{k}),\\
&\sum\limits_{j=2}^{k-1}\zeta^{(3)}(\widetilde{j},k-j)=\zeta^{(3)}(\widetilde{1},\widetilde{k-1})-\zeta^{(3)}(\widetilde{1},k-1)+\zeta^{(3)}(\widetilde{k-1},\widetilde{1})-\zeta^{(3)}(\widetilde{k-1},1)+\zeta^{(3)}(\overline{k}),\\
&\sum\limits_{j=2}^{k-1}\zeta^{(3)}(j,\overline{k-j})=\zeta^{(3)}(\overline{k-1},1)-\zeta^{(3)}(\overline{k-1},\widetilde{1})+\zeta^{(3)}(\overline{k}),\\
&\sum\limits_{j=2}^{k-1}\zeta^{(3)}(\overline{j},\overline{k-j})=\zeta^{(3)}(\overline{1},\widetilde{k-1})-\zeta^{(3)}(\overline{1},\overline{k-1})+\zeta^{(3)}(\widetilde{k-1},\overline{1})-\zeta^{(3)}(\widetilde{k-1},\widetilde{1})+\zeta(k),\\
&\sum\limits_{j=2}^{k-1}\zeta^{(3)}(\widetilde{j},\overline{k-j})=\zeta^{(3)}(\widetilde{1},k-1)-\zeta^{(3)}(\widetilde{1},\overline{k-1})+\zeta^{(3)}(\widetilde{k}),\\
&\sum\limits_{j=2}^{k-1}\zeta^{(3)}(j,\widetilde{k-j})=\zeta^{(3)}(\widetilde{k-1},1)-\zeta^{(3)}(\widetilde{k-1},\overline{1})+\zeta^{(3)}(\widetilde{k}),\\
&\sum\limits_{j=2}^{k-1}\zeta^{(3)}(\overline{j},\widetilde{k-j})=\zeta^{(3)}(\overline{1},k-1)-\zeta^{(3)}(\overline{1},\widetilde{k-1})+\zeta^{(3)}(\overline{k}),\\
&\sum\limits_{j=2}^{k-1}\zeta^{(3)}(\widetilde{j},\widetilde{k-j})=\zeta^{(3)}(\widetilde{1},\overline{k-1})-\zeta^{(3)}(\widetilde{1},\widetilde{k-1})+\zeta^{(3)}(\overline{k-1},\widetilde{1})-\zeta^{(3)}(\overline{k-1},\overline{1})+\zeta(k).
\end{align*}
\end{cor}

Recall that $\omega$ is the primitive $3$th root of unity. For positive integers $k_1,k_2$ with $k_1\geq 2$ and $a_1,a_2\in R_3$, from Lemma \ref{Lem:MZVN-MLV}, we have
\begin{align}
\zeta_{3}(k_1,k_2;a_1,a_2)=&\sum\limits_{m_1>m_2>0}\frac{\left(1+\omega^{m_1-a_1}+\omega^{2(m_1-a_1)}\right)\left(1+\omega^{m_2-a_2}+\omega^{2(m_2-a_2)}\right)}{m_1^{k_1}m_2^{k_2}}\notag\\
=&\zeta(k_1,k_2)+\omega^{-a_1}\zeta^{(3)}(\overline{k_1},k_2)+\omega^{-2a_1}\zeta^{(3)}(\widetilde{k_1},k_2)\notag\\
&+\omega^{-a_2}\zeta^{(3)}(k_1,\overline{k_2})+\omega^{-a_1-a_2}\zeta^{(3)}(\overline{k_1},\overline{k_2})+\omega^{-2a_1-a_2}\zeta^{(3)}(\widetilde{k_1},\overline{k_2})\notag\\
&+\omega^{-2a_2}\zeta^{(3)}(k_1,\widetilde{k_2})+\omega^{-a_1-2a_2}\zeta^{(3)}(\overline{k_1},\widetilde{k_2})+\omega^{-2a_1-2a_2}\zeta^{(3)}(\widetilde{k_1},\widetilde{k_2}).
\label{Eq:express-level 3}
\end{align}
Hence we can get sum formulas of multiple zeta values of level $3$ from Corollary \ref{Cor:Sum-DLV}. To state the results, we introduce some notations. For an integer $k$ with $k\geq 3$, we set
\begin{align*}
\zeta_3^{0,1}(1,\overline{k-1})=&\sum\limits_{m_1>m_2>0}\frac{(1-\omega)\omega^{m_1+2}(\omega^{m_1-2}-1)(1+\omega^{m_2-1}+\omega^{2(m_2-1)})}{m_1m_2^{k-1}},\\
\zeta_3^{1,2}(1,\overline{k-1})=&\sum\limits_{m_1>m_2>0}\frac{(1-\omega)\omega^{m_1+1}(\omega^{m_1}-1)(1+\omega^{m_2-1}+\omega^{2(m_2-1)})}{m_1m_2^{k-1}},\\
\zeta_3^{2,0}(1,\overline{k-1})=&\sum\limits_{m_1>m_2>0}\frac{(1-\omega)\omega^{m_1}(\omega^{m_1+2}-1)(1+\omega^{m_2-1}+\omega^{2(m_2-1)})}{m_1m_2^{k-1}},\\
\zeta_3^{1,0}(1,\widetilde{k-1})=&\sum\limits_{m_1>m_2>0}\frac{(1-\omega)\omega^{m_1+2}(1-\omega^{m_1-2})(1+\omega^{m_2+1}+\omega^{2(m_2+1)})}{m_1m_2^{k-1}},\\
\zeta_3^{0,2}(1,\widetilde{k-1})=&\sum\limits_{m_1>m_2>0}\frac{(1-\omega)\omega^{m_1}(1-\omega^{m_1+2})(1+\omega^{m_2+1}+\omega^{2(m_2+1)})}{m_1m_2^{k-1}},\\
\zeta_3^{2,1}(1,\widetilde{k-1})=&\sum\limits_{m_1>m_2>0}\frac{(1-\omega)\omega^{m_1+1}(1-\omega^{m_1})(1+\omega^{m_2+1}+\omega^{2(m_2+1)})}{m_1m_2^{k-1}}.
\end{align*}

Using Corollary \ref{Cor:Sum-Double of level 3} and \eqref{Eq:express-level 3}, we get the following sum formulas of double zeta values of level $3$.

\begin{cor}
For an integer $k$ with $k\geq 3$, we have
\begin{align*}
&\sum\limits_{j=2}^{k-1}\zeta_3(j,k-j)=3\zeta_3(k),\\
&\sum\limits_{j=2}^{k-1}\zeta_3(j,\overline{k-j})=\zeta_3^{2,0}(1,\overline{k-1})+\zeta_3(\overline{k-1},\widetilde{1})-\zeta_3(k-1,\widetilde{1}),\\
&\sum\limits_{j=2}^{k-1}\zeta_3(j,\widetilde{k-j})=\zeta_3^{1,0}(1,\widetilde{k-1})+\zeta_3(\widetilde{k-1},\overline{1})-\zeta_3(k-1,\overline{1}),\\
&\sum\limits_{j=2}^{k-1}\zeta_3(\overline{j},k-j)=\zeta_3(k-1,\overline{1})-\zeta_3(\overline{k-1},\overline{1}),\\
&\sum\limits_{j=2}^{k-1}\zeta_3(\overline{j},\overline{k-j})=\zeta_3^{0,1}(1,\overline{k-1}),\\
&\sum\limits_{j=2}^{k-1}\zeta_3(\overline{j},\widetilde{k-j})=\zeta_3^{2,1}(1,\widetilde{k-1})+\zeta_3(\widetilde{k-1},\widetilde{1})-\zeta_3(\overline{k-1},\widetilde{1})+3\zeta_3(\widetilde{k}),\\
&\sum\limits_{j=2}^{k-1}\zeta_3(\widetilde{j},k-j)=\zeta_3(k-1,\widetilde{1})-\zeta_3(\widetilde{k-1},\widetilde{1}),\\
&\sum\limits_{j=2}^{k-1}\zeta_3(\widetilde{j},\overline{k-j})=\zeta_3^{1,2}(1,\overline{k-1})+\zeta_3(\overline{k-1},\overline{1})-\zeta_3(\widetilde{k-1},\overline{1})+3\zeta_3(\overline{k}),\\
&\sum\limits_{j=2}^{k-1}\zeta_3(\widetilde{j},\widetilde{k-j})=\zeta_3^{0,2}(1,\widetilde{k-1}).
\end{align*}
\end{cor}

Similarly, taking all possible values of $(a_1,a_2)$ in Corollary \ref{Cor:Weightedsum-double},  we get the following weighted sum formulas of multiple $L$-values of level $3$.

\begin{cor}\label{Cor:Weightedsum-double level 3}
For an integer $k$ with $k\geq 3$, we have
\begin{align*}
&\sum\limits_{j=2}^{k-1}2^{j}\zeta(j,k-j)=(k+1)\zeta(k),\\
&\sum\limits_{j=2}^{k-1}2^{j}\zeta^{(3)}(\overline{j},k-j)=2\zeta^{(3)}(\overline{1},\widetilde{k-1})-2\zeta^{(3)}(\overline{1},k-1)
\\
&\qquad\qquad\qquad\qquad+2\zeta^{(3)}(\widetilde{k-1},\overline{1})-2\zeta^{(3)}(\widetilde{k-1},\widetilde{1})+(k-1)\zeta^{(3)}(\widetilde{k})+2\zeta(k),\\
&\sum\limits_{j=2}^{k-1}2^{j}\zeta^{(3)}(\widetilde{j},k-j)=2\zeta^{(3)}(\widetilde{1},\overline{k-1})-2\zeta^{(3)}(\widetilde{1},k-1)
+2\zeta^{(3)}(\overline{k-1},\widetilde{1})\\
&\qquad\qquad\qquad\qquad\qquad-2\zeta^{(3)}(\overline{k-1},\overline{1})+(k-1)\zeta^{(3)}(\overline{k})+2\zeta(k),\\
&\sum\limits_{j=2}^{k-1}2^{j-1}\zeta^{(3)}(j,\overline{k-j})+\sum\limits_{j=2}^{k-1}2^{j-1}\zeta^{(3)}(\overline{j},\widetilde{k-j})=\zeta^{(3)}(\overline{1},\overline{k-1})-\zeta^{(3)}(\overline{1},\widetilde{k-1})
\\
&\qquad\qquad\qquad\qquad\qquad+\zeta^{(3)}(\overline{k-1},\overline{1})-\zeta^{(3)}(\overline{k-1},\widetilde{1})+k\zeta^{(3)}(\overline{k})+\zeta^{(3)}(\widetilde{k}),\\
&\sum\limits_{j=2}^{k-1}2^{j-1}\zeta^{(3)}(j,\widetilde{k-j})+\sum\limits_{j=2}^{k-1}2^{j-1}\zeta^{(3)}(\widetilde{j},\overline{k-j})=\zeta^{(3)}(\widetilde{1},\widetilde{k-1})-\zeta^{(3)}(\widetilde{1},\overline{k-1})
\\
&\qquad\qquad\qquad\qquad\qquad+\zeta^{(3)}(\widetilde{k-1},\widetilde{1})-\zeta^{(3)}(\widetilde{k-1},\overline{1})+k\zeta^{(3)}(\widetilde{k})+\zeta^{(3)}(\overline{k}),\\
&\sum\limits_{j=2}^{k-1}2^{j-1}\zeta^{(3)}(\overline{j},\overline{k-j})+\sum\limits_{j=2}^{k-1}2^{j-1}\zeta^{(3)}(\widetilde{j},\widetilde{k-j})=\zeta^{(3)}(\overline{1},k-1)-\zeta^{(3)}(\overline{1},\overline{k-1})
\\
&\qquad\qquad\qquad\qquad+\zeta^{(3)}(\widetilde{1},k-1)-\zeta^{(3)}(\widetilde{1},\widetilde{k-1})+(k-1)\zeta(k)+\zeta^{(3)}(\overline{k})+\zeta^{(3)}(\widetilde{k}).\\
\end{align*}
\end{cor}

Using Corollary \ref{Cor:Weightedsum-double level 3} and \eqref{Eq:express-level 3}, we get the following weighted sum formulas of double zeta values of  level $3$.

\begin{cor}
For an integer $k$ with $k\geq 3$, we have
\begin{align*}
\sum\limits_{j=2}^{k-1}2^{j}\zeta_3(\overline{j},\widetilde{k-j})&=(2\omega+4)\left(\zeta^{(3)}(\overline{1},k-1)-\zeta^{(3)}(\widetilde{1},\widetilde{k-1})\right)\\
&+(2\omega-2)\left(\zeta^{(3)}(\overline{1},\overline{k-1})-\zeta^{(3)}(\widetilde{1},k-1)\right)\\
&+(4\omega+2)\left(\zeta^{(3)}(\widetilde{1},\overline{k-1})-\zeta^{(3)}(\overline{1},\widetilde{k-1})\right)\\
&+(3k-3)\zeta(k)+(3k-3)\omega\zeta^{(3)}(\overline{k})-(3k-3)(\omega+1)\zeta^{(3)}(\widetilde{k}),\\
\sum\limits_{j=2}^{k-1}2^{j}\zeta_3(\widetilde{j},\overline{k-j})&=(2\omega+4)\left(\zeta^{(3)}(\widetilde{1},k-1)-\zeta^{(3)}(\overline{1},\overline{k-1})\right)\\
&+(2\omega-2)\left(\zeta^{(3)}(\widetilde{1},\widetilde{k-1})-\zeta^{(3)}(\overline{1},k-1)\right)\\
&+(4\omega+2)\left(\zeta^{(3)}(\overline{1},\widetilde{k-1})-\zeta^{(3)}(\widetilde{1},\overline{k-1})\right)\\
&+(3k-3)\zeta(k)+(3k-3)\omega\zeta^{(3)}(\widetilde{k})-(3k-3)(\omega+1)\zeta^{(3)}(\overline{k}).
\end{align*}
\end{cor}

\end{document}